\documentclass[12 pt]{article}

\usepackage[dvips]{graphicx}
\usepackage{amsthm,amsfonts,amsmath}

\newcommand{\RR}{\mathbb{R}}
\newcommand{\ZZ}{\mathbb{Z}}
\newcommand{\kk}{\mathcal{K}}
\newcommand{\ex}{\mathbf{E}}
\newcommand{\pr}{\mathbf{P}}
\newcommand{\rr}{\stackrel {d}{=}}
\newcommand{\indic}{\mathbf{1}}

\newtheorem{thm}{Theorem}[section]

\newtheorem{prop}[thm]{Proposition}

\newtheorem{defi}[thm]{Definition}
\newtheorem{rema}[thm]{Remark}
\newtheorem{cor}[thm]{Corollary}

\allowdisplaybreaks[4]
\numberwithin{equation}{section}

\begin{document}

\begin{titlepage}
\title{\bf Diffusion local time storage}

\author{ \bf M. Kozlova\thanks{Corresponding author.
Fax: +358-2-2154865. \it E-mail address: \rm mkozlova@abo.fi}, \,
\bf P. Salminen
\vspace{0.4cm} \\
\it Department of Mathematics, \AA bo Akademi University, \\
\it Turku, FIN-20500, Finland}

\date{\today}

\maketitle

\thispagestyle{empty}

\begin{abstract}
 In this paper we study a storage process or a liquid queue
 in which the input process is the local time
 of a positively recurrent stationary diffusion in stationary state
 and the potential output takes place with a constant deterministic
 rate.
 For this storage process we find its stationary distribution and
 compute the joint distribution of the starting and
 ending times of the busy and idle periods. This work completes and
 extends to a more general setting the results in Mannersalo, Norros, and Salminen
 \cite{mannersalonorros}.
\end{abstract}
\noindent
\it MSC:
\rm  60G10; 60J55; 60J60; 60K25

\vspace{0.5cm}
\noindent
\it Keywords:
\rm
Brownian motion with drift;
Busy and idle periods;
Diffusion processes;
Palm probability;
Spectrally positive L\'{e}vy processes

\end{titlepage}

\section{Introduction} \label{intro}

Let $A=\{A_t: t \in \RR\}$ and $B=\{B_t:t \in \RR\}$ be increasing
stochastic processes. We view $A$
as the input process to a storage
or a buffer and $B$ as the potential output.
One of the standard ways to define a storage process
associated to $A$ and $B$ is via Reich's formula
 (see Reich \cite{reich58})
 $$S_t:=\sup_{-\infty<s\leq t} \{A_t-A_s-(B_t-B_s)\}, \quad t \in \RR.$$
 We refer also to Prabhu \cite{prabhu98} p. 114
 for an approach in which the storage process
 is defined as a solution of a stochastic integral equation and the input
 process $A$ is a subordinator without drift
 and $B_t=t$. In Harrison \cite{harrison85}
 p. 14 the case in which $A$ is a Brownian motion and
 $B_t= \mu t$ is studied via Skorohod's reflection equation.
Brownian motion, although not increasing,
is much analyzed input process and its appearence in these models can
be motivated in different ways, see, e.g., Harrison \cite{harrison85} p. 30,
Roberts, Mocci and Virtamo \cite{robertsmoccivirtamo96} p. 377, and
Salminen and Norros
\cite{salminennorros01}.
A new kind of a model is introduced in Mannersalo, Norros, and Salminen
\cite{mannersalonorros}, in which the local time at 0 of a reflected
Brownian motion with negative drift serves as the input process. Recall that
this local time process is increasing and continuous but the points of
its increase form a singular set.

In this paper, we consider a model of a storage or a fluid queue,
 where the input is the local time $L$ at zero of an
 arbitrary one-dimensional positively recurrent stationary diffusion and
 the potential output
 $B_t =\mu t$, $\mu>0$; thus extending the approach in
\cite{mannersalonorros}. Recall that any subordinator can be viewed as
the inverse local time of a strong Markov process. Hence, in a sense,
the local time storage models are companions of the subordinator
storage models and as such interesting objects of research.

It is shown in \cite{mannersalonorros} that the storage model studied
therein can be obtained
 via a limiting procedure from some on/off processes. Especially,
a priority discrete queueing system is constructed
which, when in heavy traffic,
  leads to a storage process with local time input.
  The proof of weak convergence in \cite{mannersalonorros} is based
   on the continuous mapping theorem and the interpretation of the Brownian local time as
  a supremum process.
 However, this approach is not applicable in our general case,
  and we do not discuss here the weak converegence aspects of our model.

 Further, in \cite{mannersalonorros} the joint distribution of
 the starting time $g_i$ and the ending time $d_i$ of
 an idle period and the marginal distributions of the
 starting time $g_b$ and the ending time $d_b$
 of a busy period are obtained.
 The \it joint \rm distribution of $g_b$ and $d_b$
 is not, however, given in \cite{mannersalonorros}. The aim of this
 paper is to
 fill this gap and to generalize the results to an \it arbitrary
 diffusion
 \rm local time storage.

 We give explicit expressions for the joint Laplace transforms of the starting
 and the ending times of the busy and idle periods.
 Surprisingly,
 the answers have a short simple form, expressed in terms of the L\'{e}vy-Khinchin
 exponent of the inverse of the local time process $L.$
 By considering a marked point process obtained from the storage process $S$
 we make the simple but crucial observation needed for finding the joint
 distributions. This observation, formulated here only for the busy
 period, is that the random variable
$(-g_b, d_b)$ is identical in law with $(UV, (1-U)V)$, where $U$ has
the uniform distribution
 on $(0,1)$ and $V$ is a positive random variable independent of $U$
 having the distribution of the length of the busy period.
 In particular, which is important here,
 we show that the joint Laplace transform of such
 variables can be easily expressed in terms of the marginal Laplace transforms.

 We remark also that in \cite{mannersalonorros} the marginal distribution of
 $-g_b$ is computed by making use of Bertoin's path decomposition
 theorem (see Bertoin \cite{bertoin91}).
 However, in this paper we  focus instead on $d_b$.
 It is seen that the formula for the first time when a spectrally positive L\'{e}vy
  process jumps over a level, see (\ref{bing}), can be applied to find
  the Laplace transform of $d_b.$

 The paper is organized as follows. In the next section we give
 some definitions
 and construct the storage process with a diffusion local time input.
 It is important to verify here that our local time process
 has stationary increments.
It is also proved that the stationary distribution of
this storage process is always
 exponential with an atom at zero - only the value of the parameter and
 the size of the atom vary when the underlying diffusion is changed.
 In Section $\ref{busyperiods}$, the Laplace transform of the ending time
 $d_b$ of a busy period observed at time zero is computed.
 In Section $\ref{idleperiods}$, we find the Laplace transform
 of the starting time $g_i$ of an
 idle period observed at time zero.
 Theorem \ref{theorem3} in Section $\ref{joint}$ can be viewed as
 the main result of the paper, where the joint distributions of
 $(-g_b,d_b)$ and $(-g_i, d_i)$ are given via their Laplace transforms.
Finally, in Section $\ref{example}$, we discuss the case treated in
\cite{mannersalonorros}. In particular, we invert the Laplace transforms to obtain the
joint densities of $(-g_i, d_i)$ and $(-g_b, d_b)$.

\section{Definitions and preliminaries} \label{preliminaries}

Let $X=\{X_t: t \geq 0\}$ be a one-dimensional time-homogeneous, regular, conservative
 diffusion living in an interval $I \subseteq \RR$. We assume, without loss
 of generality, that $0 \in I$.
The probability measure and the expectation
associated with $X$ started at $x$, are denoted by $\pr_x$ and $\ex_x$, respectively.
Recall that $X$ has a symmetric
transition density with respect to its speed
measure
$m$:
$$\pr_x (X_t \in dy)=p(t;x,y) m(dy)=p(t;y,x) m(dy).$$
We introduce also the Green function
$$G_\alpha(x,y) = \int_0^\infty e^{-\alpha t} p(t; x,y) dt, \quad x,y \in I, \: \alpha \geq 0.$$
Further, it is assumed that $X$ is positively recurrent, that is,
the speed measure is finite in $I$. Denote $M:=m\{I\} < +\infty$.
For more details on linear diffusions we refer to It\^{o} and McKean
\cite{itomckean74} and Borodin and Salminen \cite{borodinsalminen02}.

Let $X^{(1)}$ and
$X^{(2)}$ be
two copies of $X$.
We take $X_0^{(1)} = X_0^{(2)}$ with the common distribution
$$\pr(X_0^{(1)} \in dx)= \pr(X_0^{(2)} \in dx)= \frac{m(dx)}{M}=:\widehat m(dx)$$
but otherwise we let $X^{(1)}$ and
$X^{(2)}$ be independent. We remark that $\widehat m$ is the
stationary probability distribution of $X.$
Define $Z=\{Z_t: t \in \RR\}$ as
$$Z_t:=\left \{
    \begin{array}{cc}
       X_{-t}^{(1)} & \mbox{if $t \leq 0$}, \\
       X_t^{(2)} & \mbox{if $t \geq 0$}.
     \end{array}
     \right.$$
The process $Z$ thus defined is a stationary process in stationary state.

Let $\{L_t^{(1)} : t \geq 0\}$ and $\{L_t^{(2)} : t \geq 0\}$
be the local times at zero of the processes $X^{(1)}$ and $X^{(2)}$, respectively. We
choose the normalization of $L^{(i)}$, $i=1,2$, with respect to the speed measure, that is,
    \begin{equation}
         L_t^{(i)}:= \lim_{\varepsilon \searrow 0}
         \frac{1}{m\{(-\varepsilon, \varepsilon)\}}
         \int_0^t \indic_{(- \varepsilon, \varepsilon)}
         (X_s^{(i)}) \,ds,\quad i=1,2.
         \label{eqloc}
    \end{equation}
\noindent
Using $L^{(1)}$ and $L^{(2)}$, we define the local time process
$L=\{ L_t: t \in \RR \}$ via
$$L_t:=\left \{
    \begin{array}{cc}
       -L_{-t}^{(1)} & \mbox{if $t \leq 0$}, \\
       L_t^{(2)} & \mbox{if $t \geq 0$}.
     \end{array}
     \right.$$
The process $L$ is an increasing continuous process passing through zero at time zero.
Next we give an important property of $L$, which allows us
to construct a stationary storage process.

\begin{prop} \label{stat}
The process $L$ has stationary increments, i.e., for all $t>~s$,
$$L_t-L_s \rr L_{t-s}.$$
Moreover,
\begin{equation}\label{mean}
\ex (L_t) = t /M.
\end{equation}
\end{prop}
\noindent
\textit{Proof:} Consider first the case $t>s>0$. Then, using that $Z$ is a stationary
process in stationary state, we have
\begin{eqnarray*}
  \pr(L_t-L_s \in dl) & = & \int_I \widehat m(dx) \pr(L_t-L_s \in dl \, | \, Z_s=x) \\
  & = & \int_I \widehat m(dx) \pr_x(L_{t-s} \in dl) \\
  & = & \pr(L_{t-s} \in dl).
\end{eqnarray*}
The case $s < t < 0$ is treated similarly.
Assume next that $s<0<t$. Then by the conditional independence of $X^{(1)}$ and
$X^{(2)}$ given $Z_0$ we can write
\begin{eqnarray*}
  \ex(e^{-\alpha (L_t-L_s)}) & = & \ex (e^{-\alpha(L_t^{(1)} +L_{-s}^{(2)})}) \\
  & = & \int_I \widehat m(dx) \ex_x (e^{-\alpha L_t^{(1)}})
  \ex_x (e^{-\alpha L_{-s}^{(2)}}).
\end{eqnarray*}
On the other hand, consider
\begin{eqnarray*}
  \ex(e^{-\alpha L_{t-s}}) & = & \ex(e^{\alpha (L_{t-s} -L_{-s} +
  L_{-s})}) \\
  & = & \int_I \widehat m(dx) \int_I \ex_x(e^{-\alpha L_{-s}} \, ; \, Z_{-s} \in dy)
  \ex(e^{\alpha (L_{t-s} -L_{-s})} \, | \, Z_{-s} = y).
\end{eqnarray*}
Let now $\tau$ be an exponentially distributed random
variable with parameter $\alpha$
 independent of $X^{(2)}$, and define the
killed process
$$X_t^{\bullet}:=\left \{
    \begin{array}{cc}
       X_{t}^{(2)} & \mbox{if $L_t^{(2)} < \tau$}, \\
       \triangle & \mbox{if $L_t^{(2)} \geq \tau$},
     \end{array}
     \right.$$
 where $\triangle$ is a cemetery point.
Then (see Borodin and Salminen \cite{borodinsalminen02} p. 28 or It\^{o} and McKean
\cite{itomckean74} p. 179-183) $X_t^{\bullet}$
is a diffusion having the same
speed measure as $X^{(2)}$. Let $p^{\bullet}(t; x,y)$  be its
(symmetric) transition density with
respect to $m$. We clearly have
\begin{eqnarray*}
   \pr_x (X_t^{\bullet} \in dy) & = & \pr_x (X_t^{(2)} \in dy, L_t^{(2)}< \tau)\\
   & = & \ex_x (e^{-\alpha L_t^{(2)}} \,; \,X_t^{(2)} \in dy)\\
   & = & p^{\bullet}(t; x,y)\,m(dy).
\end{eqnarray*}
From the symmetry of $p^{\bullet}(t; x,y)$ it follows that
$$\int_I m(dx) \ex_x(e^{-\alpha L_{-s}}\, ; \, Z_{-s} \in dy)
= m(dy) \ex_y (e^{-\alpha L_{-s}}),$$
and, therefore,
\begin{eqnarray*}
  \ex(e^{-\alpha L_{t-s}}) & = & \int_I \widehat m(dy) \ex_y(e^{-\alpha L_{-s}})
  \ex(e^{-\alpha(L_{t-s} - L_{-s})} \,|\, Z_{-s} = y) \\
  & = & \int_I \widehat m(dy) \ex_y(e^{-\alpha L_{-s}})
   \ex_y (e^{-\alpha L_t})\\
  & = & \ex (e^{-\alpha(L_t-L_s)}),
\end{eqnarray*}
as claimed.
Finally, to show $(\ref{mean})$ recall the formula (see It\^{o}
and McKean \cite{itomckean74} p. 175 or Borodin and Salminen
\cite{borodinsalminen02} p. 21)
$$\ex_x(L_t) = \int_0^t p(s; x,0) \, ds, \quad t >0.$$
Consequently,
\begin{eqnarray*}
 \ex(L_t) & = & \int_I \widehat m(dx) \int_0^t ds \, p(s;x,0)\\
 & = & \int_0^t ds \int_I \widehat m(dx)\, p(s;x,0) \\
 & = & \int_0^t ds \int_I \widehat m(dx)\, p(s;0,x) \\
 & = & t / M,
\end{eqnarray*}
because
$$\int_I m(dx)\, p(s;0,x)=1.$$ \qed

\begin{defi} \label{defstorl}
\rm      The process $S=\{S_t: t \in \RR\}$ defined as
      $$S_t:=\sup_{-\infty < s \leq t}\{L_t-L_s- \mu (t-s)\}, \quad t \in \RR$$
      is called a \it storage process with local time input, constant service rate $\mu$, and
      unbounded buffer associated with the process $Z$.
\end{defi}

Because the process $L$ has stationary increments it follows from Definition
$\ref{defstorl}$ that $S$ is a stationary process (in stationary state). Our next
task is to find its stationary distribution, that is, the distribution of $S_0$.
Notice that $S_0$ is determined by the process $X^{(1)}$. With a slight and short
abuse of our notations, let
$\{L^{(1)}_t: t \geq 0\}$ be the local time of $X^{(1)}$ started at zero
normalized as in $(\ref{eqloc})$. Then the process $A^{(1)}$ given by
$$A^{(1)}_t:=\inf\{s: L^{(1)}_s > t\}, \quad t \geq 0$$
is a subordinator (a right-continuous increasing
process on $\RR_+$ with independent and stationary increments).
When normalizing as in (\ref{eqloc}) we have
(see It\^{o} and McKean \cite{itomckean74}, p. 214) that
\begin{equation} \label{levy}
   \ex_0 (\exp \{-\alpha A^{(1)}_t\}) = \exp \left\{ -\frac{t}{G_\alpha(0,0)} \right\},
\end{equation}
and, therefore,
\begin{equation*}
   \ex_0 \left(\exp\{-\alpha (A^{(1)}_t-t/ \mu) \}\right) =
   \exp\left \{-t\left(\frac{1}{G_\alpha(0,0)}-\frac{\alpha}{\mu} \right) \right\}
   = \exp \{t \psi(\alpha) \},
   \end{equation*}
where the function
\begin{equation} \label{levy1}
\psi(\alpha)=\frac{\alpha}{\mu} - \frac{1}{G_\alpha(0,0)}
\end{equation}
is called the L\'{e}vy-Khintchin exponent of
the spectrally positive L\'{e}vy process $\{A^{(1)}_t - t/ \mu \,: \, t \geq 0\}$.

\begin{prop} \label{rasprs0}
 Assume that $\mu > 1/M$, in other words $\mu > \ex (L_1)$. Then the equation
 \begin{equation} \label{psi}
  \psi(\alpha)=0
 \end{equation}
has  a unique positive solution $\alpha^*$. Moreover, \vspace{0.2cm}\\
 \bf 1) \it $\displaystyle{\pr(S_0>t| Z_0 =x)= \ex \left( e^{-\alpha^* H_0}
 | Z_0=x\right)
   e^{-\frac{\alpha^*}{\mu} t}, \quad t \geq 0, \; x \in I}$,
   \vspace{0.2cm}\\
   where $H_0:=\inf \{t \geq 0: \, Z_t = 0\}$;
   \vspace{0.2cm}\\
   \bf 2) \it the process $S$ is stationary with the stationary
    distribution given by
    $$\pr(S_t>y) = \frac{1}{M \mu} e^{-\frac{\alpha^*}{\mu} y},\quad y \geq 0,$$
    and $ \pr(S_t = 0) = 1 - \frac{1}{M \mu}$.
\end{prop}
\noindent
\textit{Proof:}
The stationarity follows from Proposition $\ref{stat}$.
The other claims follow from the fact that
the spectrally positive L\'{e}vy process $\{A^{(1)}_t - t/ \mu \,: \, t \geq
0\}$ drifts to $+\infty$ and its infimum is exponentially distributed in
$\RR_-$ with parameter $\alpha^*.$ For further details, see
\cite{salminen93}. \qed

\begin{rema} \label{rema1}
  \rm
\bf 1) \rm The Laplace transform of $H_0$ appearing above can be
  expressed as
     \begin{equation} \label{green}
         \ex_x \left( e^{-\alpha H_0}\right)=
         \frac{G_\alpha(x,0)}{G_\alpha(0,0)},
       \end{equation}
    cf. Borodin and Salminen \cite{borodinsalminen02}, p. 18. \\
   \bf 2) \rm Recall (cf. Bingham \cite{bingham75} p. 720) that for $\alpha \geq \alpha^*$, the function $\psi(\alpha)$
   is positive and increasing and, hence, has an inverse which
   we denote by $\alpha \mapsto
   \eta(\alpha)$. \\
   \bf 3) \rm From Proposition $\ref{rasprs0}$ it is clear that
$S_0$, conditionally on $S_0>0$, is exponentially distributed with parameter
$\alpha^* / \mu$.
\end{rema}
 We assume throughout the paper that $\mu > 1/M$.
Next we define the subject of our main interest in this work, the on-going
busy and idle periods of $S$.

\begin{defi} \label{defidle}
\rm    If $S_0 = 0$ then the random variables
    $$g_i = \sup\{t<0: S_t>0\} \quad \mbox{\rm and} \quad
    d_i= \inf\{t>0: S_t>0\}$$
    are called the \it starting time \rm and the \it ending time, \rm respectively,
    of the \it on-going idle period \rm at time zero.
    If $S_0 > 0$ then the random variables
    $$g_b = \sup\{t<0: S_t=0\} \quad \mbox{\rm and} \quad
    d_b= \inf\{t>0: S_t=0\}$$
    are called the \it starting time \rm  and the \it ending time, \rm  respectively,
    of the \it on-going busy period \rm at time zero.
\end{defi}

\section{Busy periods} \label{busyperiods}

In this section we compute the Laplace transform of the ending time of the on-going
busy period. It is shown in Mannersalo, Norros, and Salminen
\cite{mannersalonorros}
using a point process view on the busy and idle periods
that $-g_b$ and $d_b$ have
the same law. In the present paper
(in Section $\ref{joint}$) we develop this approach in order to find the joint distribution
of $-g_b$ and $d_b$.

From Proposition $\ref{rasprs0}$
it is seen that the probability that there is a busy period at time zero is $1/(M \mu)$.

\begin{thm} \label{theorem1}
Let $\eta$ be the inverse of $\psi$ (cf. Remark $\ref{rema1}.2$). Then
    \begin{equation} \label{formula1}
      \ex(e^{-\alpha d_b} \,|\, S_0>0)=
      \frac{\alpha^*}{\eta(\alpha/ \mu)}.
     \end{equation}
\end{thm}

\begin{figure}[!h]
\begin{center}
\includegraphics[width=12cm,height=8cm,angle=0]{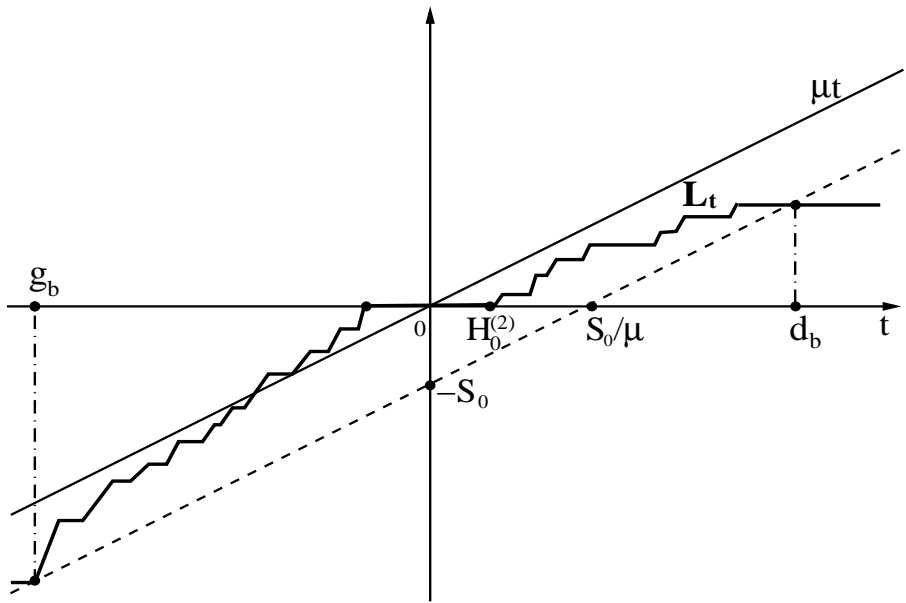}
\caption{The on-going busy period at time zero.}\label{ris1}
\end{center}
\end{figure}

The rest of this section is devoted to the proof of Theorem $\ref{theorem1}$.
Assuming that there is a busy period at 0 we have (cf. \cite{mannersalonorros})
$$
S_t=L_t-\mu t + S_0,\quad g_b\leq t\leq d_b.
$$
Let
$H_0^{(2)}$ be the first hitting time of zero for the process
$X^{(2)}$.
Because
$$S_0 = \sup_{s \leq 0} \{\mu s -L_s\}=
\sup_{s \leq 0} \{\mu s +L_{-s}^{(1)}\}$$
it follows that $S_0$
and $H_0^{(2)}$ are conditionally independent, given $Z_0$ .
The local time $L_t=0$ when $0 \leq t \leq H_0^{(2)}$.
Hence,
  $$S_t = S_0 - \mu t,\quad 0 \leq t  \leq \min\{H_0^{(2)}, S_0/ \mu\}.
$$
To proceed consider, therefore, two cases:

\begin{description}
  \item[\bf 1)] $S_0 / \mu < H_0^{(2)}$,
  \item[\bf 2)] $S_0 / \mu > H_0^{(2)}$ (see Figure $\ref{ris1}$).
\end{description}
In the case \bf 1) \rm we obviously have
\begin{equation} \label{dbs}
  \{ S_0 / \mu < H_0^{(2)}\, , \, S_0>0\} = \{d_b=S_0/ \mu \, , \, S_0>0\}.
\end{equation}
To analyze the case
\bf 2)\rm, we study in detail some relationships between $S_0$ and
$H_0^{(2)}$.
Firstly, introduce
\begin{equation} \label{j}
J(\alpha , \beta): = \int_I m(dx) \ex_x(e^{-\alpha H_0})
\ex_x(e^{-\beta H_0}), \quad \alpha \geq 0,\:  \beta\geq 0.
\end{equation}

\begin{prop} \label{prop1}
 \begin{equation} \label{shh}
    \ex \left(e^{-\alpha(S_0 / \mu - H_0^{(2)}) -\beta H_0^{(2)}} \, , \, S_0 / \mu>H_0^{(2)} \right) =
     \frac{\alpha^*}{\alpha+ \alpha^*}
 \frac{1}{M} J(\alpha^* , \alpha^* +\beta).
  \end{equation}
\end{prop}
\noindent
\textit{Proof}: Using Proposition $\ref{rasprs0}$ and the
conditional independence of $S_0$ and
 $H_0^{(2)}$ given $Z_0$, we compute
 the joint Laplace transform as follows:
 \begin{eqnarray*}
   & & \ex \left(e^{-\alpha(S_0 / \mu -H_0^{(2)}) -\beta H_0^{(2)}}
   \, , \, S_0 / \mu>H_0^{(2)} \right) \\
& & \hspace{0.7cm} = \int_I \widehat m(dx)
\ex_x \left(e^{-\alpha (S_0 / \mu-H_0^{(2)})-\beta H_0^{(2)}} \, ,
\, S_0 / \mu>H_0^{(2)} \right) \\
 & & \hspace{0.7cm} = \int_I \widehat m(dx) \int_0^\infty \ex_x
 \left(e^{-\alpha (S_0 / \mu-t)} e^{-\beta t} \, , \, S_0 / \mu>t \right)
    \pr_x (H_0^{(2)} \in dt) \\
& & \hspace{0.7cm} = \int_I \widehat m(dx) \int_0^\infty
 e^{(\alpha -\beta)t}\pr_x (H_0^{(2)} \in dt)
    \int_t^\infty e^{-\alpha u}\pr_x(S_0 / \mu \in du) \\
 & & \hspace{0.7cm} = \int_I \widehat m(dx) \int_0^\infty
 e^{(\alpha -\beta)t}\pr_x (H_0^{(2)} \in dt)
\int_t^\infty e^{-\alpha u} \alpha^* e^{- \alpha^* u} \ex_x(e^{-\alpha^* H_0}) du \\
& & \hspace{0.7cm} = \frac{\alpha^*}{\alpha+\alpha^*} \frac{1}{M} \int_I m(dx)
 \int_0^\infty
 e^{(\alpha -\beta)t}\pr_x (H_0^{(2)} \in dt)
 \ex_x(e^{-\alpha^* H_0}) e^{-(\alpha+\alpha^*)t}\\
   & & \hspace{0.7cm} = \frac{\alpha^*}{\alpha+\alpha^*}  \frac{1}{M}\int_I m(dx)
   \ex_x(e^{-\alpha^* H_0}) \ex_x(e^{-(\alpha^*+\beta)H_0}),
   \end{eqnarray*}
as claimed. \qed

An immediate consequence of Proposition $\ref{prop1}$ is the following
\begin{cor} \label{cor1}
 \bf 1) \rm The random variables $S_0 / \mu-H_0^{(2)}$
 and $H_0^{(2)}$ are conditionally independent,
 given that $S_0 / \mu >H_0^{(2)}$. \\
 \bf 2) \rm  The conditional distribution of $S_0 / \mu-H_0^{(2)}$, given that $S_0 / \mu>H_0^{(2)}$, is
       exponential with parameter $\alpha^*$.  \\
 \bf 3) \rm  $\displaystyle{\ex \left(e^{-\alpha H_0^{(2)}} \, ,  \, S_0 / \mu > H_0^{(2)} \right) =
       \frac{J(\alpha^*, \alpha^* +\alpha)}{M}}.$
       \vspace{0.2cm}
         \\
\bf 4) \rm
$\displaystyle{\ex \left( e^{-\alpha S_0 / \mu} \, , \, S_0 / \mu>H_0^{(2)} \right) =
        \frac{\alpha^*}{\alpha+\alpha^*}  \frac{J(\alpha^*, \alpha^* +
        \alpha)}{M}}.$
\end{cor}

\begin{prop} \label{prop2}
The functions $J(\alpha^*, \alpha^*+\alpha)$  and $\psi(\alpha)$
defined in $(\ref{j})$ and $(\ref{levy1})$, respectively,  satisfy the
relationship:
  \begin{equation} \label{I}
    J(\alpha^*, \alpha^* +\alpha) = \frac{1}{\mu} -\frac{\psi(\alpha+\alpha^*)}{\alpha}.
  \end{equation}
\end{prop}

\textit{Proof}: Using formula $(\ref{green})$ and the Chapman-Kolmogorov equation,
we compute
for any positive $\alpha_1 \ne \alpha_2$,
\begin{eqnarray*}
  & & \int_I m(dx) \, \ex_x (e^{-\alpha_1 H_0}) \ex_x (e^{-\alpha_2 H_0})\\
  & & \hspace{0.2cm} =
      \int_I  m(dx) \, G_{\alpha_1} (x,0) G_{\alpha_2} (x,0) \, /(G_{\alpha_1} (0,0) G_{\alpha_2} (0,0))\\
  & & \hspace{0.2cm} =
      \int_I m(dx) \int_0^\infty dt \, e^{-\alpha_1 t} p(t;x,0)
       \int_0^\infty ds \, e^{-\alpha_2 s} p(s;x,0)
       /(G_{\alpha_1} (0,0) G_{\alpha_2} (0,0)) \\
& & \hspace{0.2cm} =
      \int_0^\infty dt\int_0^\infty ds \, e^{-\alpha_1 t- \alpha_2 s}
      \int_I m(dx) p(t;x,0) p(s;x,0)  /(G_{\alpha_1} (0,0) G_{\alpha_2}(0,0)) \\
& & \hspace{0.2cm} =
      \int_0^\infty dt \int_0^\infty ds \,e^{-\alpha_1 t- \alpha_2 s}
       p(t+s,0,0)
   \,   /(G_{\alpha_1} (0,0) G_{\alpha_2}(0,0)) \\
  & & \hspace{0.2cm} =
       \int_0^\infty du\, \left( \int_0^u e^{(\alpha_1-\alpha_2) s} ds \right) e^{-\alpha_1 u}
       p(u; 0,0)\,  /(G_{\alpha_1} (0,0) G_{\alpha_2}(0,0)) \\
  & & \hspace{0.2cm} =
       \frac {1}{\alpha_1-\alpha_2} \int_0^\infty (e^{-\alpha_2 u} - e^{-\alpha_1 u}) p(u;0,0) du
       \, /(G_{\alpha_1} (0,0) G_{\alpha_2}(0,0)) \\
   & & \hspace{0.2cm} = \frac{1}{\alpha_1-\alpha_2} \left( \frac{1}{G_{\alpha_1} (0,0)}
   -\frac{1}{G_{\alpha_2}(0,0)} \right).
\end{eqnarray*}
Since
$$\frac{1}{G_\alpha(0,0)}= \frac{\alpha}{\mu} - \psi(\alpha),$$
we have
\begin{equation} \label{ipsi}
 J(\alpha_1, \alpha_2)=
 \frac{1}{\mu}-\frac{\psi(\alpha_1) -
   \psi(\alpha_2)}{\alpha_1-\alpha_2}.
\end{equation}
Substituting $\alpha^*$ for $\alpha_1$ and $\alpha+\alpha^*$ for
$\alpha_2$ in $(\ref{ipsi})$ and using that
$\psi(\alpha^*)=0$ give $(\ref{I})$. \qed

To proceed with the proof of Theorem $\ref{theorem1}$ write
\begin{eqnarray}
  \ex(e^{-\alpha d_b} \, |  \, S_0>0)& = &
  \ex(e^{-\alpha d_b} \, , \, S_0 / \mu < H_0^{(2)} \,|\, S_0>0)
  \nonumber \\
   & & \hspace{0.2cm} + \,  \ex(e^{-\alpha d_b} \, ,
   \, S_0 / \mu > H_0^{(2)}\, |\, S_0>0). \label{dbsumma}
\end{eqnarray}
For the first term we obtain by $(\ref{dbs})$, Corollary $\ref{cor1}.4$, and
Proposition $\ref{rasprs0}$
\begin{eqnarray}  \label{lap1}
  & & \ex(e^{-\alpha d_b} \, ,  \, S_0 / \mu  <  H_0^{(2)} \, | \, S_0>0) \nonumber \\
  & &  \hspace{0.7cm} = \ex(e^{-\alpha S_0 / \mu} \, , \,  S_0 / \mu  <  H_0^{(2)} \, |  \, S_0>0) \nonumber \\
  & & \hspace{0.7cm} = \ex(e^{-\alpha S_0 / \mu} \, |  \, S_0 > 0)-
  \ex(e^{-\alpha S_0 / \mu} \, , \, S_0 / \mu  >  H_0^{(2)} \, | \, S_0>0) \nonumber \\
  & & \hspace{0.7cm} = \ex(e^{-\alpha S_0 / \mu} \,| \, S_0 > 0)-
  \frac{\ex (e^{-\alpha S_0 / \mu} \, , \, S_0 / \mu  >  H_0^{(2)})}{\pr(S_0>0)}  \nonumber \\
  & & \hspace{0.7cm} = \frac{\alpha^*}{\alpha + \alpha^*}
   \big( 1- \mu J(\alpha^*, \alpha^* + \alpha)\big).
\end{eqnarray}
To compute the second term in $(\ref{dbsumma})$ consider the process
$$\widetilde{X}= \{X_{H_0^{(2)} + t}: t \geq 0\}.$$
By the strong Markov property $\widetilde{X}$  has the
same law as $X^{(2)}$ started at zero. Let
$\widetilde{L}= \{\widetilde{L_t}: t \geq 0\}$ be the local time at zero of $\widetilde{X}$,
that is, the local time of $X^{(2)}$ started at zero.
Let $\widetilde{A}$
be the right continuous inverse of $\widetilde{L}$ and introduce
$$T^+_t:=\inf\{s: \widetilde{A}_s-s/ \mu > t\}, \quad t>0.$$
Observe next that in the case $S_0 / \mu > H_0^{(2)}$ (see Figure
 \ref{ris1}) we have
$$d_b= \frac{1}{\mu}\, T^+_{S_0 / \mu - H_0^{(2)}}+ \frac{S_0}{\mu}.$$
Using this representation and conditional independence
stated in Corollary
 $\ref{cor1}.1$, it is seen that
\begin{eqnarray*}
  && \ex \left(e^{-\alpha d_b} \, , \, S_0 / \mu > H_0^{(2)} \, \big{|} \, S_0>0 \right) \\
  & & \hspace{0.7cm} = \ex \left(e^{-\alpha d_b} \, \big{|}
   \, S_0 / \mu > H_0^{(2)} \right)
   \pr \left(S_0 / \mu > H_0^{(2)} \, \big{|} \, S_0 >0 \right) \\
  & & \hspace{0.7cm} = \ex \left (\exp \left\{-\alpha
  \left(\frac{1}{\mu}\, T^+_{S_0 / \mu-H_0^{(2)}}+ \frac{S_0}{\mu} \right)\right\}
   \, \Big{|} \, S_0 / \mu > H_0^{(2)}  \right) \\
  & & \hspace{1.7cm} \times
  \pr \left(S_0 / \mu > H_0^{(2)} \, \big{|} \, S_0 >0 \right)\\
  & & \hspace{0.7cm} = \ex \left( \exp\left\{ -\alpha \left(
  \frac{1}{\mu} \,
   T^+_{S_0 / \mu-H_0^{(2)}}+ \frac{S_0}{\mu}-H_0^{(2)} \right) \right\}
  \,\Big{|}\, S_0 / \mu > H_0^{(2)} \right) \\
  & & \hspace{1.7cm} \times \ex \left(e^{-\alpha H_0^{(2)}} \, \big{|} \,
   S_0 / \mu > H_0^{(2)} \right)
  \pr \left(S_0 / \mu > H_0^{(2)} \,\big{|}\, S_0 >0 \right) \\
   & & \hspace{0.7cm} = \ex \left(\exp\left\{-\alpha  \left(
   \frac{1}{\mu} \,
   T^+_{S_0 / \mu-H_0^{(2)}}+ \frac{S_0}{\mu} -H_0^{(2)} \right)\right\}
  \, \Big{|} \, S_0 / \mu > H_0^{(2)} \right) \\
    & & \hspace{1.7cm} \times \ex \left(e^{-\alpha H_0^{(2)}} \, , \,
    S_0 / \mu > H_0^{(2)} \right)
    \big{/} \,\pr \left(S_0>0 \right).
\end{eqnarray*}
From Corollary $\ref{cor1}.2$ and the fact that $T^+$ is independent of $S_0$
and $H_0^{(2)}$, we obtain
\begin{eqnarray*}
&&\ex \left(\exp\left\{-\alpha \left( \frac{1}{\mu}T^+_{S_0 / \mu-H_0^{(2)}}
+S_0 / \mu-H_0^{(2)} \right) \right\} \, \Big{|} \, S_0 / \mu > H_0^{(2)} \right) \\
& & \hspace{0.7cm} = \int_0^\infty
\ex \left(\exp\left\{-\frac{\alpha}{\mu} T^+_t  - \alpha t \right\} \right)
 \pr \left(S_0 / \mu-H_0^{(2)} \in dt \,\big{|}\, S_0 / \mu > H_0^{(2)} \right) \\
& &  \hspace{0.7cm} = \alpha^* \int_0^\infty
\ex \left(\exp\left\{-\frac{\alpha}{\mu} T^+_t \right\}\right)
e^{-(\alpha+ \alpha^*)t}\, dt.
\end{eqnarray*}
To compute the last integral we recall the following
formula for the double Laplace transform of the first time when the spectrally positive L\'{e}vy process
$\{\widetilde{A}_s -s/ \mu \, : \, s \geq 0\}$
jumps over the level $t>0$
(see Bingham \cite{bingham75}, p. 732):
\begin{equation} \label{bing}
   \int_0^\infty \ex \left(e^{-\gamma T^+_t } \right) e^{-\beta t} dt =
   \frac {1}{\beta} \left( 1- \frac{\gamma \,(\eta(\gamma) - \beta)}
   {\eta(\gamma)(\gamma - \psi(\beta))} \right),
\end{equation}
where $\alpha \mapsto \psi(\alpha)$ is given by $(\ref{levy1})$ and
$\alpha \mapsto \eta(\alpha)$ is the inverse of
$\psi(\alpha),\: \alpha \geq \alpha^*$.
Combining this with the results of Proposition $\ref{rasprs0}$ and
Corollary $\ref{cor1}$ yields
\begin{eqnarray} \label{formula2}
  & & \ex \left(e^{-\alpha d_b} \, , \, S_0/ \mu > H_0^{(2)} \, \big{|}\, S_0>0\right) \nonumber \\
  & & \hspace{0.5cm} = \frac{\alpha^*}{\alpha+\alpha^*}
  \left( 1- \frac{\frac{\alpha}{\mu}
  \left(\eta(\frac{\alpha}{\mu})-(\alpha+\alpha^*) \right)}{\eta(\frac{\alpha}{\mu})
   \left(\frac{\alpha}{\mu} -
  \psi(\alpha+\alpha^*) \right)} \right) \mu \, J(\alpha^*,\alpha^* +\alpha).
\end{eqnarray}
Finally, recalling $(\ref{dbsumma})$,
putting $(\ref{formula2})$ together with $(\ref{lap1})$, and using $(\ref{I})$
complete the proof of Theorem $\ref{theorem1}$.
\pagebreak

\pagebreak

\section{Idle periods} \label{idleperiods}

\begin{figure}[!h]
\begin{center}
\includegraphics[width=12cm,height=8cm,angle=0]{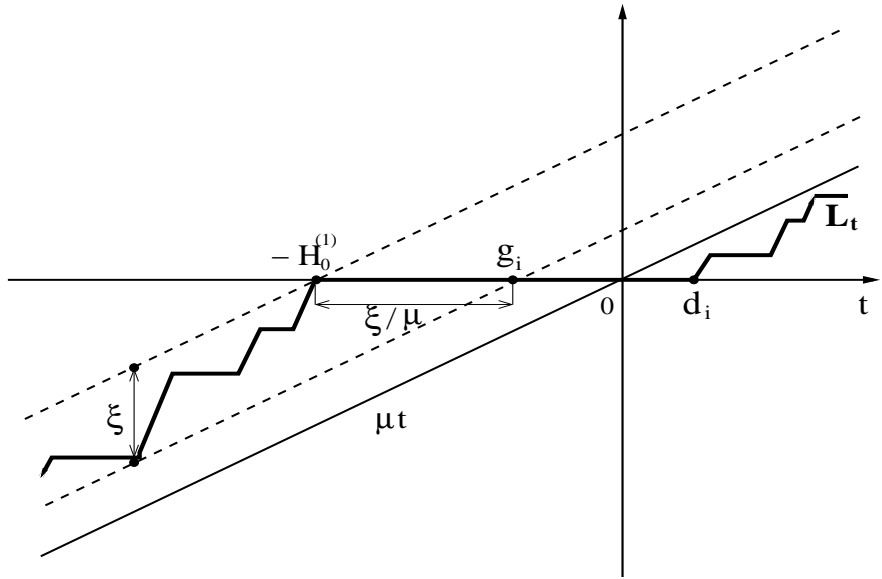}
\caption{The on-going idle period at time zero.}\label{ris2}
\end{center}
\end{figure}

Let $H_0^{(1)}$ be the first hitting time of zero for the process
$X^{(1)}$ and
introduce
$$\xi:=\sup_{t \leq 0} \{ \mu t - L_{-H_0^{(1)}+t}\}.$$
By the strong Markov property of $X^{(1)}$ and Proposition
 $\ref{rasprs0}.1$, the
random variable
$\xi$ is exponentially distributed
with parameter $\alpha^* / \mu$ (cf. Remark $\ref{rema1}.3$) and independent of $H_0^{(1)}$.
The local time $L_t$ is equal to zero when $-H_0^{(1)} \leq t \leq 0$,
and
for these values of $t$
$$S_t = \sup_{s \leq 0} \{\mu s - L_{s+t}\}.$$
Since $\xi=S_{-H_0^{(1)}} >0$ it is clear that $-H_0^{(1)} \leq g_i \leq 0$,
as displayed in Figure $\ref{ris2}$.
Further, we have
(see Figure $\ref{ris2}$) that
$$\{S_0=0\} = \{\xi / \mu < H_0^{(1)}\},$$
and  that the conditional distribution of $-g_i$, given that $S_0=0$, is the same as the
conditional distribution of $H_0^{(1)}- \xi / \mu$, given that $\xi / \mu <H_0^{(1)}$.

\begin{thm} \label{theorem2}
  For $\alpha \ne \alpha^*$,
  \begin{equation} \label{theorem2f}
  \ex (e^{\alpha g_i} \,|\, S_0=0) = \frac{1}{\alpha}
  \frac{\mu \alpha^* \psi(\alpha)}{(M \mu -1) (\alpha-\alpha^*)}.
  \end{equation}
\end{thm}
\noindent
\textit{Proof:}        Using the independence of $\xi / \mu$ and $H_0^{(1)}$ and
$(\ref{ipsi})$, we compute
\begin{eqnarray*}
  & & \ex (e^{\alpha g_i} \, , \, S_0=0) \\
  & & \hspace{0.5cm} = \ex(e^{-\alpha (H_0^{(1)} - \xi / \mu)} \, , \,
   \xi / \mu < H_0^{(1)}) \\
  & & \hspace{0.5cm} = \int_I \widehat m(dx) \,
  \ex_x(e^{-\alpha (H_0^{(1)} - \xi / \mu)} \, , \, \xi / \mu < H_0^{(1)}) \\
   & & \hspace{0.5cm} =  \int_I \widehat m(dx) \int_0^\infty \pr_x(H_0^{(1)} \in dt)
   e^{-\alpha t} \ex_0 (e^{\alpha \xi / \mu}\, , \, \xi / \mu <t) \\
   & & \hspace{0.5cm} = \int_I  \widehat m(dx) \int_0^\infty \pr_x(H_0^{(1)} \in dt)
   e^{-\alpha t} \frac{\alpha^*}{\alpha^* - \alpha} (1-e^{-(\alpha^* - \alpha)t})
   \end{eqnarray*}
  \begin{eqnarray*}
   & & \hspace{0.5cm} = \frac{\alpha^*}{M(\alpha^*-\alpha)} \int_I m(dx)
   \left( \ex_x(e^{-\alpha H_0}) - \ex_x(e^{-\alpha^* H_0})\right) \\
   & & \hspace{0.5cm} = \frac{\alpha^*}{M(\alpha^*-\alpha)}
   \left( \frac{1}{\mu} -\frac{\psi(\alpha)}{\alpha} -\frac{1}{\mu} +
   \frac{\psi(\alpha^*)}{\alpha^*}\right) \\
    & & \hspace{0.5cm} = \frac{\alpha^* \psi(\alpha)}{M \alpha(\alpha-\alpha^*)}.
\end{eqnarray*}
Noting that
$$\pr(S_0=0) = 1- \frac{1}{M \mu}$$
gives $(\ref{theorem2f})$. \qed
\begin{rema} \label{remark2}
  \rm  In case $\alpha = \alpha^*$, $(\ref{theorem2f})$
  should be read as
  $$\ex (e^{\alpha^* g_i} \,|\, S_0=0) = \frac{\mu}{M \mu -1} \psi'(\alpha^*).$$
\end{rema}
\pagebreak

\pagebreak

\section{The joint distribution of the starting and the ending times of
 the on-going busy and idle periods} \label{joint}

In this section, as our main result, we compute the joint distribution of the starting time and the ending
time of the busy and idle periods. We remark that it is possible to compute
the joint distribution directly in the case of the idle period, but we have not been able
to do this in the busy period case.
To resolve this difficulty we use
the theory of Palm probability and
properties of a special class of two-dimensional random variables to show that it is
enough to know only the Laplace transforms of $d_b$ (respectively, $-g_i$) in order to
obtain the joint Laplace
transforms of $d_b$ and $-g_b$ (respectively, $d_i$ and $-g_i$).

\begin{thm}\label{theorem3}
\bf 1) \rm For $\alpha \ne \beta$,
      \begin{equation} \label{jointbusy}
      \ex \left( e^{\alpha g_b-\beta d_b} \,|\, S_0>0 \right)
       = \frac{\alpha^*}{\alpha-\beta} \left( \frac{\alpha}{\eta(\alpha / \mu)}
       -\frac{\beta}{\eta(\beta / \mu)}\right).
      \end{equation}
\bf 2) \rm For $\alpha \ne \beta$,
     \begin{equation} \label{jointidle}
      \ex \left( e^{\alpha g_i-\beta d_i} \,|\, S_0=0 \right)
      = \frac{\alpha^* \mu}{(M \mu-1) (\alpha-\beta)}
      \left( \frac{\psi(\alpha)}{\alpha-\alpha^*} -
      \frac{\psi(\beta)}{\beta-\alpha^*}\right).
     \end{equation}
\end{thm}

\begin{rema} \label{remark3}
\rm  Letting $\beta \to \alpha$ in $(\ref{jointbusy})$ and
$(\ref{jointidle})$ gives us
  the Laplace transforms of the lengths of the busy
  and idle periods
  $$\ex \left( e^{-\alpha (d_b-g_b)} \,|\, S_0>0 \right)
  = \alpha^* \frac{d}{d \alpha} \left( \frac{\alpha}{\eta(\alpha / \mu)}\right)$$
  and
 $$\ex \left( e^{-\alpha (d_i-g_i)} \,|\, S_0>0 \right)=
  \frac{\alpha^* \mu}{M \mu-1} \frac{d}{d \alpha} \left(
   \frac{\psi(\alpha)}{\alpha-\alpha^*}\right),$$
respectively.
\end{rema}

To prove
Theorem $\ref{theorem3}$ we
use the point process approach to the storage process $S$ proposed in Mannersalo, Norros,
and Salminen \cite{mannersalonorros}.
When $S$ hits zero it stays there for
a positive amount of time
 a.s. Hence, we can construct from $S$  a stationary marked
point process $N=\{(T_n, Z_n): n \in \ZZ\}$ with a finite
intensity $\lambda$, where for every $n$, $T_n$ and $T_{n+1}$ are
the starting and ending times, respectively, of an idle or a busy period and
$Z_n$ is the mark associated with point $T_n$. We let $Z_n=0$
if $T_n$ is a starting
point for an idle period, otherwise $Z_n = 1$.
Let $\pr_N$ and $\mathcal{F}$ denote the probability measure
and the natural $\sigma$-algebra, respectively, induced by $N$.
We assume that $\{T_n: n \in \ZZ\}$ are ordered so that $-\infty <
   \cdots < T_{-1} < T_0 \leq 0 < T_1 <T_2 < \cdots + \infty$.
Let $\theta_t$ be the shift operator and  $\pr_N^0$  the Palm probability,
 which is related to
the probability measure $\pr_N$  via
Ryll-Nardzewski and Slivnyak's formula:
\begin{equation}
   \pr_N(A) = \lambda \int_0^\infty \pr_N^0(T_1 >t, \theta_t \in A) \, dt,
   \quad A \in \mathcal{F}.
    \label{eqinverse}
\end{equation}
Recall that the Palm probability $\pr_N$ can be interpreted as the conditional
probability, given that there is a point at time zero, i.e.,
   \begin{equation} \label{eqp}
    \pr_N^0(\{T_0 = 0\}) = 1.
   \end{equation}
For more details about these concepts, see
Baccelli and Br\'{e}maud
\cite{baccellibremaud03} or Mannersalo, Norros, and Salminen
\cite{mannersalonorros}.
From relationship (\ref{eqinverse}) we obtain now the following property
of the joint distribution of $T_0$ and $T_1$.

\begin{prop} \label{proptt}
For $z=0$ or $1$ and $v \geq 0, \: w \geq 0$,
\begin{equation}
         \pr(T_1>v, -T_0>w, \,Z_0 =z) = \lambda \int_{v+w}^\infty
         \pr_N^0(T_1>t, \,Z_0 =z)\,dt.
         \label{eqtt1}
  \end{equation}
\end{prop}
\noindent
\textit{Proof:} To establish formula $(\ref{eqtt1})$
we only need to take $A=\{T_1>v, -T_0>w, \, Z_0=z\}$ in (\ref{eqinverse})
and use that for all $t \in [0, T_1)$,
$$Z_0 \circ \theta_t = Z_0,$$
$$T_1 \circ \theta_t=T_1-t,$$ and
 $$-T_0 \circ \theta_t=t$$
$\pr_N^0$-a.s. (for $T_0$, see ($\ref{eqp})$). \qed

Consider now the on-going busy period of $S$. There are obvious analogous formulas
for the idle period.
By  (\ref{eqtt1}), we obtain
\begin{eqnarray}
    \pr(-g_b>v, d_b>w\, , \, S_0>0 ) & = &
    \pr(T_1>v, -T_0>w\, ,\, Z_0=1) \nonumber \\
    & = & \lambda \int_{v+w}^\infty \pr_N^0(T_1>u\, , \, Z_0=1)\,du. \label{eqidlepalm}
\end{eqnarray}
Consequently, letting
$$F_b^c(v,w):= \pr(-g_b>v, d_b>w\,|\, S_0>0 )$$
it is seen that
\begin{equation} \label{fk}
 F_b^c(v,w) = F_b^c(v+w,0).
\end{equation}
Two-dimensional positive random variables with the property
$(\ref{fk})$ form a special class of random
variables studied, e.g.,  in Salminen and Vallois \cite{salminenvallois03}.
The alternative definition of this class, denoted by $\kk$,  is
as follows.

\begin{defi}\label{defk}  \rm A two-dimensional positive random variable
$(X,Y)$ belongs to the \it class $\kk$  \rm if
$(X,Y)$ has the same distribution as $(UV,(1-U)V)$, where the random variable
$U$ has the uniform distribution on $(0,1)$, and $V$ is an arbitrary positive
random variable independent of $U$.
\end{defi}

\begin{rema} \label{schur}
  \rm In reliability theory the functions with the property
   $(\ref{fk})$
   are called Schur-constant functions
   (see, e.g., Bassan and
  Spizzichino \cite{bassan}).
\end{rema}

The next proposition describes the structure of the joint density of an element in $\kk$
(for the proof, see
Salminen and Vallois \cite{salminenvallois03}).
\begin{prop}\label{propk}
Let $(X,Y)\in \kk$  be such that there exists the density $F'_V(v)=:p(v)$.  Then
$$\pr (X \in dx, Y \in dy) = f(x+y) \,dx \,dy,$$
where
$$f(v)=v^{-1}p(v), \quad v >0.$$
\end{prop}

We have also the following characterization of elements in the class $\kk$ in terms
of the Laplace transforms.

\begin{prop}\label{propkl}
Let $X$ and $Y$ be positive random variables. Then $(X,Y) \in \kk$ if and only if there
exists a positive random variable $V$ such that for all $\alpha \ne \beta$,
   \begin{equation} \label{formulakl}
   \ex \left(e^{-\alpha X-\beta Y} \right) = \frac{1}{\alpha-\beta} \int_\beta^\alpha
   \ex \left(e^{-\gamma V} \right) d \gamma.
   \end{equation}
In particular, $V \rr X+Y$.
\end{prop}
\noindent
\textit{Proof:} $\Rightarrow)$ Let $(X,Y) \in \kk$. Then there
exists a positive $V$ such that
$$(X,Y) \rr (UV, (1-U)V),$$
where $U \sim U(0,1)$ is
 independent of $V$.
Taking $\alpha > \beta$ we compute
\begin{eqnarray*}
    \ex \left(e^{ -\alpha X - \beta Y} \right) & = &
    \ex  \left(e^{ -\alpha U V-\beta (1-U) V} \right)\\
    & = & \ex \left( e^{-(\alpha-\beta) U V-\beta V} \right)\\
    & = & \int_0^{\infty} \ex \left(e^{-(\alpha-\beta) U v-\beta v} \right) F_V(dv)\\
    & = & \frac {1}{\alpha - \beta}
    \int_0^\infty \frac{e^{-\beta v} - e^{-\alpha v}}{v} F_V(dv)\\
    & = &\frac{1}{\alpha - \beta} \int_{0}^\infty \left(
    \int_{\beta}^{\alpha} e^{-v \gamma} d \gamma \right) F_V(dv) \\
     & = & \frac{1}{\alpha - \beta} \int_{\beta}^{\alpha}
     \ex \left(e^{- \gamma V} \right) d \gamma,
\end{eqnarray*}
as claimed in $(\ref{formulakl})$.   \\
$\Leftarrow)$ Now suppose that we have a pair $(X,Y)$ of random variables such that
$(\ref{formulakl})$  holds.
Letting $\alpha \to \beta$ in $(\ref{formulakl})$ and using the continuity of
the Laplace transforms,
we get that $X+Y$ has the same law as $V$.
Further, let $U \sim U(0,1)$ be independent of $V$.
Then $(UV, (1-U)V) \in \kk$ and by
the sufficiency part of the proof,
 for $\alpha \ne \beta$,
  $$\ex \left(e^{-\alpha UV-\beta (1-U)V} \right) = \frac{1}{\alpha-\beta} \int_\beta^\alpha
   \ex \left(e^{-tV} \right) dt. $$
Hence, by uniqueness of the Laplace transforms
$$(X,Y) \rr (UV, (1-U)V) \in \kk.$$ This completes the proof.
\qed

\begin{rema} \label{corclassk}
 \rm If $(X,Y) \in \kk$ and the Laplace transform of $X$ is known,
we can easily compute the joint Laplace transform. Indeed, denote
$I(\alpha):= \int_0^\alpha \ex \left( e^{- \gamma V}\right) d \gamma$. Then
setting $\beta=0$ in formula $(\ref{formulakl})$, we get
$$\ex \left(e^{-\alpha X} \right) = \frac{1}{\alpha} I(\alpha)$$
and hence,
 \begin{eqnarray*}
   \ex \left(e^{-\alpha X-\beta Y} \right) & = & \frac{1}{\alpha-\beta}
   \int_\beta^\alpha  \ex \left( e^{- \gamma V}\right) d \gamma \\
   & = & \frac{1}{\alpha-\beta} \left( I(\alpha) - I(\beta) \right).
  \end{eqnarray*}
\end{rema}

We conclude now the proof of Theorem $\ref{theorem3}$. By $(\ref{eqtt1})$, the two-dimensional random variable $(-g_b, d_b) \in \kk$
and the marginal Laplace
transform of $d_b$ is given by $(\ref{formula1})$. Hence we can use Remark $\ref{corclassk}$ to get the
joint distribution of $-g_b$ and $d_b$.
From Theorem $\ref{theorem1}$ we have
$$\ex(e^{-\alpha d_b} \, | \, S_0>0) = \frac{\alpha^*}{\eta(\alpha/ \mu)}.$$
Therefore, the joint Laplace transform of $-g_b$ and $d_b$ is
$$\ex(e^{\alpha g_b -\beta d_b} \,|\, S_0>0) =
\frac{\alpha^*}{\alpha - \beta} \left(\frac{\alpha}{\eta(\alpha/ \mu)}
-\frac{\beta}{\eta(\beta / \mu)} \right).$$
The corresponding formula for $(-g_i, d_i)$ is obtained similarly. The proof
of Theorem $\ref{theorem3}$ is now complete.

\section{Example: reflected Brownian motion with negative drift}   \label{example}
Let $Z = \{Z_t: t \in \RR\}$ be a reflected Brownian motion with drift $-c<0$  in stationary state
living on $I=[0, \infty)$. Let
$$S_t:=\sup_{s \leq t}\{L_t-L_s-(t-s)\}$$
 be the local time storage process (as introduced in
 Section $\ref{preliminaries}$), associated with $Z$.
In this case we have $$m(dx)= 2 e^{-2 c x} dx,$$ and $M:=m\{I\} = 1/ c$.
The storage process $S$ is well-defined if and only if $0<c<1$.
 For a reflected Brownian motion with drift $-c$ the Green function
at $(0,0)$ is
$$G_\alpha(0,0)= \frac{1}{\sqrt{2 \alpha+ c^2}-c},$$
and the function $\psi$ (cf. (\ref{levy1})) takes the form
$$\psi(\alpha)= \alpha-\sqrt{2 \alpha+c^2} + c, \quad \alpha \geq 0$$
Hence,
$$\alpha^* = 2(1-c)$$ and
$$\eta(\alpha)= \alpha+1-c+\sqrt{2\alpha+(1-c)^2}.$$
Using $(\ref{jointbusy})$, we compute
\begin{eqnarray} \label{busybrown}
 &&\ex \left(e^{\alpha g_b -\beta d_b}\,|\, S_0>0 \right) \nonumber \\
 & & \hspace{0.2cm} = \frac{2(1-c)}{\alpha-\beta}
 \left( \frac{\alpha}{\alpha+1-c+\sqrt{2\alpha+(1-c)^2}} -
 \frac{\beta}{\beta+1-c+\sqrt{2\beta+(1-c)^2}}\right) \nonumber \\
 & & \hspace{0.2cm} = \frac{1}{\alpha-\beta}
 \left( \frac{4(1-c)}{\sqrt{2 \beta +(1-c)^2}
 +1+ c} - \frac{4(1-c)}{\sqrt{2 \alpha +(1-c)^2}
 +1+ c} \right )   \nonumber\\
 & & \hspace{0.2cm} = \frac{8(1-c)}{\sqrt{2 \alpha+(1-c)^2}+ \sqrt{2 \beta+(1-c)^2}}
 \nonumber  \\
 & &  \hspace{1.5cm} \times
 \frac{1}{(\sqrt{2 \alpha+(1-c)^2}+1+c)(\sqrt{2 \beta+(1-c)^2}+1+c)} \nonumber \\
 & & \hspace{0.2cm} =:F(\alpha, \beta; 1-c).
\end{eqnarray}
\noindent
For the on-going idle period at zero, using $(\ref{jointidle})$,
 we obtain after some cancellations that
$$\ex \left(e^{\alpha g_i -\beta d_i}\,|\, S_0>0 \right) =
F(\alpha, \beta; c).$$
This formula can also be found in \cite{mannersalonorros}
but $(\ref{busybrown})$ is a new result.

To find the joint density of $(-g_b, d_b)$, given that $S_0>0$, first find the density of
$V:=-g_b+d_b$, given that $S_0>0$. Setting $\beta = \alpha$ in the right-hand side
of $(\ref{busybrown})$,
 we get
\begin{equation} \label{lapv}
  \ex \left( e^{-\alpha (d_b -g_b)} \,|\, S_0>0 \right) =
  \frac{4 (1-c)}{\sqrt{2 \alpha +(1-c)^2} (\sqrt{2 \alpha +(1-c)^2}+1+c)^2}.
\end{equation}
\noindent
Taking the inverse Laplace transform of $(\ref{lapv})$ (cf.
Erdèlyi \cite{erdelyi54}, p. 234)
we obtain the density of the length of the busy period $d_b - g_b$, given that $S_0>0$, as
$$p_c (v) = 2(1-c) \sqrt{\frac{2 v}{\pi}} e^{-\frac{(1-c)^2}{2} v}
- 4 (1-c^2) v e^{2 c v} \Phi \left( -\frac{1+c}{2}\sqrt{v} \right), \; v>0,$$
where $\Phi(v)$ is the standard normal distribution function.
Note that the density of the length of the idle period $d_i - g_i$, given that $S_0=0$, is $p_{1-c} (v), \, v >0$.
Using Proposition $\ref{propk}$, we get the joint density of the starting and
the ending points of the busy period as
  \begin{eqnarray*}
  && \pr(-g_b \in dx, d_b \in dy \,|\, S_0>0)   \\
  & & \hspace{0.7cm} =  2 (1-c) \sqrt{\frac{2}{\pi (x+y)}}
   e^{-\frac{(1-c)^2}{2} (x+y)}  \\
   && \hspace{1.2cm} -4(1-c^2) e^{2 c (x+y)} \Phi
     \left(-\frac{1+c}{2} \sqrt{x+y}\right)\, dx \, dy.
\end{eqnarray*}

\paragraph{Acknowledgements}

The authors thank Ilkka Norros for a short but stimulating
discussion, and  Bruno Bassan and Fabio Spizzichino for
pointing out the connection of the class $\kk$ and
the Schur-constant functions.

\bibliography{art}
\bibliographystyle{plain}

\end{document}